\documentclass[preprint,3p,12pt]{elsarticle}

\usepackage{lineno,hyperref}
\modulolinenumbers[5]

\usepackage{hyperref}

\usepackage[all]{xy}

\usepackage{epstopdf}

\usepackage{mathrsfs}
\usepackage{amssymb}
\usepackage{amsfonts}
\usepackage{latexsym}
\usepackage{amsmath}
\usepackage{xcolor}
\usepackage{wrapfig}
\usepackage{floatflt}

\usepackage{graphicx}

\def\lb{\label}

\newcommand{\er}[1]{\textrm{(\ref{#1})}}

\newtheorem{theorem}{\bf Theorem}[section]

\newtheorem{definition}[theorem]{\bf Definition}

  \def\cA{{\mathcal A}}

\def\G{\Gamma}         
         
\def\D{\Delta}  \def\cF{{\mathcal F}}

    \def\cI{{\mathcal I}}

\def\l{\lambda}

\def\s{\sigma}  \def\cR{{\mathcal R}}

\def\O{\Omega}

\def\Z{{\mathbb Z}}    \def\R{{\mathbb R}}   \def\C{{\mathbb C}}

\def\lt{\biggl}                  \def\rt{\biggr}
\def\ol{\overline}               \def\wt{\widetilde}

\let\ge\geqslant                 \let\le\leqslant

\def\sm{\setminus}               
\def\ss{\subset}                 \def\ts{\times}
\def\pa{\partial}                
                 \def\ev{\equiv}
        
\def\el2{\ell^{\,2}}             \def\1{1\!\!1}

\def\det{\mathop{\mathrm{det}}\nolimits}

\def\dim{\mathop{\mathrm{dim}}\nolimits}

\let\ge\geqslant
\let\le\leqslant

\newcommand{\ca}{\begin{cases}}
\newcommand{\ac}{\end{cases}}
\newcommand{\ma}{\begin{pmatrix}}
\newcommand{\am}{\end{pmatrix}}
\def\eq{\begin{equation}}
\def\qe{\end{equation}}
\def\[{\begin{equation}}
\def\]{\end{equation}}

\begin{document}

\begin{frontmatter}

\title{Application of matrix-valued integral continued fractions to spectral problems on periodic graphs}
\date{\today}

\author
{Anton A. Kutsenko}

\address{Department of Mathematics, Aarhus University, Aarhus, DK-8000,
Denmark; email: akucenko@gmail.com}

\begin{abstract}
We show that spectral problems for periodic operators on lattices
with embedded defects of lower dimensions can be solved with the
help of matrix-valued integral continued fractions. While these
continued fractions are usual in the approximation theory, they are
less known in the context of spectral problems. We show that the
spectral points can be expressed as zeroes of determinants of the
continued fractions. They are also useful in the study of inverse
problems (one-to-one correspondence between spectral data and
defects). Finally, the explicit formula for the resolvent in terms
of the continued fractions is also provided. We apply some of our
results to the Schr\"odinger operator acting on the graphene with
line and point defects.
\end{abstract}


\end{frontmatter}


\section{Introduction}
\setcounter{equation}{0}

\begin{figure}[h]
\begin{minipage}[h]{0.49\linewidth}
\center{\includegraphics[width=0.8\linewidth]{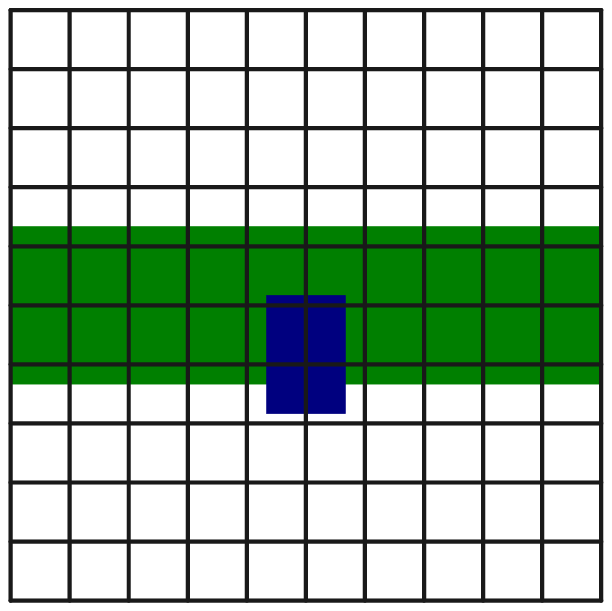}} \\ a)
\end{minipage}
\hfill
\begin{minipage}[h]{0.49\linewidth}
\center{\includegraphics[width=0.8\linewidth]{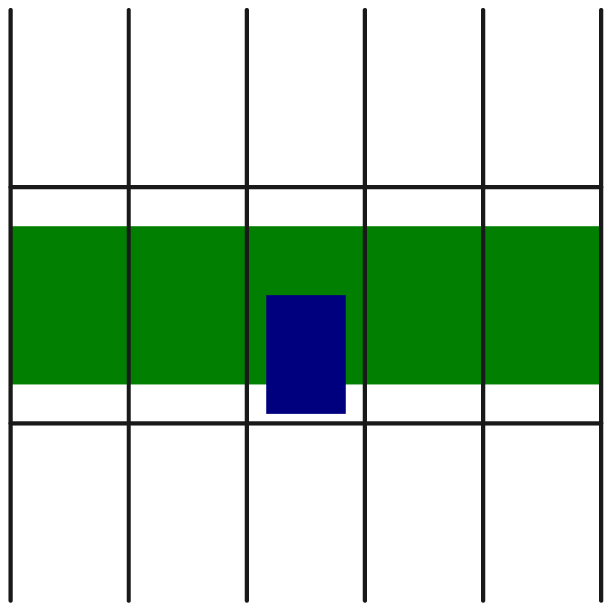}} \\ b)
\end{minipage}
\caption{2D periodic lattice with 1D periodic defect sublattice
(green strip) and local (0D) defect (blue rectangle). The net shows
different choices of  unite cells: new unite cells in (b) consist of
eight initial unite cells shown in (a).} \label{fig1}
\end{figure}

Discrete periodic operators with defects have different applications
in various problems of physics and mechanics, see discussions in
\cite{K1,K3} and an example below. The spectrum of purely periodic
operators acting on purely periodic lattices corresponds to the
extended eigenfunctions which have no attenuation along any
direction in the lattice. If we add a defect sublattice of lower
dimension to the purely periodic lattice then we obtain new spectral
components which correspond to the eigenfunctions bounded along the
defect and exponentially decreasing in the perpendicular directions.
For example, free surfaces of periodic structures and various
wave-guides embedded into the periodic lattices are defects of lower
dimensions. The corresponding extended eigenfunctions are called as
surface and guided modes. They are of high interest in various
problems of the propagation of radio (light) and acoustic waves. In
the present paper we show how matrix-valued integral continued
fractions can be used for determining the spectrum, the resolvent
and other characteristics of periodic operators with defects. We
apply these results to the Schr\"odinger operator acting on the
graphene with line and point defects.

It is shown in \cite{Kjmaa} that discrete periodic operators acting
on discrete $N$-dimensional periodic graphs (lattices) with embedded
periodic subgraphs (defect sublattices) of smaller dimensions
$N-1,...,0$ are unitarily equivalent to integral operators of a
special form. The corresponding unitary map is called a
Fourier-Floquet-Bloch transformation (FFB). In the simplest case of
a single point unite cell ($M=1$), it is the Fourier series which
allows us to replace infinite discrete sequences with functions of
continuous variables. Let
\[\lb{002}
 L^2_{N,M}\ev L^2([0,1]^N,\C^M)
\]
be the Hilbert space of square-integrable vector-valued (if $M>1$)
functions defined on $[0,1]^N$. As it is shown in \cite{Kjmaa},
after applying FFB the periodic operators with defects take the form
\[\lb{001}
 \cA:L^2_{N,M}\to L^2_{N,M},\ \ \cA={\bf A}_0\cdot+{\bf A}_1\langle{\bf
 B}_1\cdot\rangle_{1}+...+{\bf A}_N\langle{\bf
 B}_N\cdot\rangle_{1,N},
\]
where ${\bf A}\ev{\bf A}({\bf k})$, ${\bf B}\ev{\bf B}({\bf k})$ are
continuous matrix-valued functions on ${\bf k}=(k_i)\in[0,1]^N$ of
sizes
\[\lb{003}
 \dim({\bf A}_0)=M\ts M,\ \ \dim({\bf B}_j)=M_j\ts M,\ \ \dim({\bf
 A}_j)=M\ts M_j,\ \ j\ge1
\]
with some positive integers $M_j$; the dot $\cdot$ in \er{001}
denotes the place of operator arguments ${\bf u}\in L^2_{N,M}$; the
integrals $\langle\cdot\rangle_{i,j}$ in \er{001} are defined as
follows
\[\lb{e001b}
 \langle\cdot\rangle_{i,j}=\int\limits_{[0,1]^{j-i+1}}\cdot
 dk_idk_{i+1}...dk_j\ (j>i),\ \ \
 \langle\cdot\rangle_{i}\ev\langle\cdot\rangle_{i,i}=\int_{0}^1\cdot dk_i.
\]
Note that for purely periodic lattices without defects, the periodic
operators are unitarily equivalent to $\cA={\bf A}_0\cdot$, see,
e.g., \cite{KK1}. The presence of defects leads to additional
integral terms ${\bf A}_j\langle{\bf B}_j\cdot\rangle_{1,j}$, see
\cite{Kjmaa}. The matrix-valued functions ${\bf A}$, ${\bf B}$
\er{001} depend on the structure of the periodic lattice and its
defect sublattices. Roughly speaking, ${\bf A}_j$, ${\bf B}_j$ show
how the defect sublattice of the dimension $N-j$ is embedded into
the substrate lattice of the dimension $N$ and how it is related to
other defect sublattices. The periodicity of the lattice means that
there is some unit cell which can be periodically translated to
cover our lattice. The value $M$ \er{003} is exactly the number of
nodes inside the unite cell. While we have different unite cells for
different sublattices it is supposed that they have the same periods
and hence we can choose the common unite cells. The choice of the
unit cell is not unique. In particular, we can increase the size of
the unite cell twice or of integer times, see Fig.
\ref{fig1}.(a),(b). The new obtained cell can also be considered as
a unite cell for our periodic lattice. This procedure increases the
size $M$ of the space $L^2_{N,M}$ but, at the same time, it can
simplify the structure of the operator \er{001}. Namely, if the
defect cells are fully integrated into the new unit cells then the
operator \er{001} takes the form
\[\lb{004}
 \cA={\bf A}_0\cdot+{\bf A}_1\langle\cdot\rangle_1+...+{\bf
 A}_N\langle\cdot\rangle_{1,N},
\]
where $M\ts M$ matrix-valued functions ${\bf A}$ satisfy
\[\lb{005}
 {\bf A}_0\ev{\bf A}_0({\bf k}_0),\ \ {\bf k}_0=(k_1,...,k_N);\ \ \  {\bf
 A}_j\ev{\bf A}_j({\bf k}_j),\ \ {\bf k}_j=(k_{j+1},...,k_N).
\]
In particular, ${\bf A}_N$ is a constant matrix. Roughly speaking,
the independence of ${\bf A}_j$ on $k_1,...,k_j$ means that we can
cover the defect sublattice of the dimension $N-j$ by the shifted
unite cells, where we shift one "zero" unite cell along the defect
sublattice and no shifts in perpendicular directions (FFB of
perpendicular shifts are $e^{2\pi ik_r}\cdot$, $r\le j$, see
\cite{Kjmaa}) are used. For example, the unite cells of spring-mass
models considered in \cite{K1} satisfy these requirements.

Operators \er{004}-\er{005} form a linear subspace in the space
(algebra, see \cite{Kjmaa}) of all periodic operators with parallel
defects \er{001}-\er{003}. For the operators \er{004}-\er{005} the
procedure of finding the spectrum (obtained in \cite{Kjmaa} for the
general case) can be refined. Note that the Hermitian adjoint
operator of $\cA$ \er{004} has the same form
\[\lb{005a}
 \cA^*={\bf A}_0^*\cdot+\langle{\bf A}_1^*\cdot\rangle_1+...+\langle{\bf
 A}_N^*\cdot\rangle_{1,N}={\bf A}_0^*\cdot+{\bf A}_1^*\langle\cdot\rangle_1+...+{\bf
 A}_N^*\langle\cdot\rangle_{1,N},
\]
since ${\bf A}$ satisfy \er{005}.

\begin{definition}\lb{D1} Let $\l\in\C$ and let matrix-valued
functions ${\bf A}$ be of the form \er{005}. Define the
matrix-valued continued fractions ${\bf F}_j\ev{\bf F}_j(\l,{\bf
k}_j)$ by
\[\lb{e010}
 {\bf F}_0={\bf A}_0-\l{\bf I},\ \ \ {\bf F}_1={\bf
 A}_1+\left\langle\frac{\bf I}{{\bf A}_0-\l{\bf
 I}}\right\rangle_1^{-1},
\ \ \
 {\bf F}_2={\bf A}_2+\left\langle\frac{\bf I}{{\bf
 A}_1+\left\langle\frac{\bf I}{{\bf A}_0-\l{\bf
 I}}\right\rangle_1^{-1}}\right\rangle_2^{-1}
\]
and so on ${\bf F}_{j}={\bf A}_j+\langle{\bf
F}_{j-1}^{-1}\rangle_j^{-1}$ up to $j=N$. The matrix ${\bf I}$ is
the identity matrix.
\end{definition}

The next theorem gives us the expression of the spectrum of the
operator $\cA$ \er{004} in terms of the continued fractions
\er{e010}.

\begin{theorem}\lb{T1}
Let an operator $\cA$ be of the form \er{004}-\er{005}. Then the
spectrum of $\cA$ is
\[\lb{006}
 \s(\cA)=\bigcup_{j=0}^N\s_j,\ \ \s_j=\{\l:\det{\bf G}_j(\l,{\bf k}_j)=0\ for\ some\ {\bf
 k}_j\in[0,1]^{N-j}\},
\]
where
\[\lb{007}
 {\bf G}_0\ev{\bf F}_0,\ \ {\bf G}_j\ev\langle{\bf
 F}_{j-1}^{-1}\rangle_j{\bf F}_j={\bf I}+\langle{\bf
 F}_{j-1}^{-1}\rangle_j{\bf A}_j,\ \ j\ge1.
\]
\end{theorem}

{\bf Remark on the computation of the spectrum $\s(\cA)$.} On the
zero step we determine $\s_0$. On the first step, it is convenient
to define ${\bf G}_1$ for $\l\in\C\sm\s_0$ because it is not well
defined for $\l\in\s_0$ (the matrix ${\bf F}_0$ is non-invertible,
see \er{006}-\er{007}). Zeroes of ${\bf G}_1$ determine the spectral
component $\s_1$ which is disjoint from $\s_0$. And so on, on the
$j$-th step we can define ${\bf G}_j$ for
$\l\in\C\sm(\s_0\cup...\cup\s_{j-1})$. Zeroes of ${\bf G}_j$
determine $\s_j$ which is disjoint from all $\s_{j-1},...,\s_0$. In
general, it is also possible to define ${\bf G}_j$ for
$\l\in\s_0\cup...\cup\s_{j-1}$ but it leads to non-disjoint sets
$\s_j$ (of course, the total spectrum $\s=\cup\s_j$ does not
change). It is not important here but we use some of these arguments
below.

Let us briefly discuss some aspects of inverse spectral problems.
Due to \er{006} the set of matrix-valued functions $\{{\bf
G}_j\}_{j=0}^N$ can be considered as spectral data for our operator
$\cA$ \er{004}. We just need to show that $\{{\bf G}_j\}_{j=0}^N$
determine the matrix-valued functions $\{{\bf A}_j\}_{j=0}^N$
uniquely.

\begin{theorem}\lb{T2} The following identities hold true
\[\lb{008}
 {\bf A}_0=\l{\bf I}+{\bf G}_0,\ \ {\bf A}_j=\left\langle({\bf G}_0...{\bf
 G}_{j-1})^{-1}\right\rangle_{1,j}^{-1}({\bf G}_j-{\bf I}),\ j\ge1,
\]
and
\[\lb{009}
 {\bf A}_0=\l{\bf I}+{\bf F}_0,\ \ {\bf A}_j={\bf F}_{j}-\langle{\bf
F}_{j-1}^{-1}\rangle_j^{-1},\ j\ge1.
\]
\end{theorem}
Identities \er{e010}, \er{007}-\er{009} show that there are
one-to-one mappings between the sets of $M\ts M$ matrix-valued
functions
\[\lb{010}
\xymatrix{\{{\bf F}_j\}_{j=0}^N\ar@{<->}[rd]\ar@{<->}[rr]&&{\{{\bf
G}_j\}_{j=0}^N}\ar@{<->}[ld]\\&\{{\bf A}_j\}_{j=0}^N}.
\]
Thus, $\{{\bf F}_j\}_{j=0}^N$ and $\{{\bf G}_j\}_{j=0}^N$ can be
considered as spectral data in inverse spectral problems for
operators $\cA$ \er{004}. In particular, {\it if we have a set of
$M\ts M$ matrix-valued functions $\{{\bf G}_j(\l,{\bf
k}_j)\}_{j=0}^N$. Then this set corresponds to some operator $\cA$
\er{004}-\er{005} if and only if the matrix-valued functions ${\bf
A}_j$ \er{008} do not depend on $\l$. In this case, the obtained
${\bf A}_j$ are the components of the operator $\cA$ and its
spectrum satisfies \er{006}.}

Note also that the set $\{{\bf F}_j\}_{j=0}^N$ has some advantages
over the set $\{{\bf G}_j\}_{j=0}^N$: if all ${\bf A}_j$ are
self-adjoint then $\cA$ \er{004} is self-adjoint and all ${\bf F}_j$
are also self-adjoint for real $\l$. At the same time, ${\bf G}_j$
are usually non-self-adjoint.

Next remark can be useful in applications: It is not difficult to
see that instead of ${\bf G}_j$ in \er{006} we can also use
$\ol{{\bf G}}_j$ defined by
\[\lb{007a}
 \ol{{\bf G}}_0\ev{\bf F}_0,\ \ \ol{{\bf G}}_j\ev{\bf F}_j\langle{\bf
 F}_{j-1}^{-1}\rangle_j={\bf I}+{\bf A}_j\langle{\bf
 F}_{j-1}^{-1}\rangle_j,\ \ j\ge1.
\]
Thus the matrices $\ol{{\bf G}}_j$ allow us to define the spectrum
as well as the matrices ${\bf G}_j$. The analogues of Theorem
\ref{T2} and \er{010} are also fulfilled for $\ol{{\bf G}}_j$.

{\bf Remark on the resolvent.} Let $\cI$ be the identity operator.
Let $\cA$ be of the form \er{004}-\er{005} and let
$\l\not\in\s(\cA)$. The following explicit formula for the inverse
operator (resolvent) holds true
\[\lb{011}
 \cR(\l)\ev(\cA-\l\cI)^{-1}={\bf D}_0\cdot-{\bf H}_1{\bf A}_1\langle{\bf
 D}_1\cdot\rangle_1-...-{\bf H}_N{\bf A}_N\langle{\bf
 D}_N\cdot\rangle_{1,N},
\]
where the matrix-valued functions ${\bf H}$, ${\bf D}$ are defined
by
\[\lb{012}
 {\bf H}_{j+1}=(\ol{\bf G}_j...\ol{\bf G}_0)^{-1},\ \ {\bf D}_j=({\bf G}_0...{\bf
 G}_j)^{-1}.
\]
Recall that the matrix-valued functions ${\bf G}_j$ (and $\ol{{\bf
G}}_j$) determine the spectrum of $\cA$, see \er{006}-\er{007} (and
\er{007a}). If all ${\bf A}_j^*={\bf A}_j$ are self-adjoint and
$\l\in\R\sm\s(\cA)$ then $\cA$ and its resolvent are self-adjoint,
and
\[\lb{013}
 \cR(\l)={\bf H}_1\cdot-\sum_{r=1}^N{\bf
 H}_r{\bf A}_r\langle{\bf H}_{r+1}^*\cdot\rangle_{1,r}={\bf
 D}_0\cdot-\sum_{r=1}^N{\bf D}_{r-1}^*{\bf A}_r\langle{\bf
 D}_r\cdot\rangle_{1,r}.
\]
Let us provide some examples where the resolvent plays important
role. In a functional calculus of operators: we can express an
analytic function $f$ (defined on some domain $\O\ss\C$) of the
operator $\cA$ as the Cauchy's integral $f(\cA)=(2\pi
i)^{-1}\oint_{\pa\O}f(\l)\cR(\l)d\l$. In a problem with energy
sources: we can explicitly express the wave-function ${\bf u}$
corresponding to the source ${\bf f}$ as ${\bf u}=\cR(\l){\bf f}$,
where $\l$ depends on the energy of the source. We have used such
explicit expressions for the defect detection, see \cite{K3}.

{\bf The case $M=1$.} In this special case, all ${ A}$, ${ F}$, ${
G}$, ${ H}$, ${ D}$ are scalar functions (for scalars we do not use
bold fonts). Assume also that $A_j$ are real functions. Then the
operator $\cA$ \er{004}-\er{005} is self-adjoint. Define $\l_0({\bf
k}_0)\ev A_0({\bf k}_0)$, ${\bf k}_0\in[0,1]^N$. It is a standard
function but it is more convenient to look at it as a set. By
induction, suppose that for $r=0,...,j-1$ we have already defined
the sets $\l_r\ev\l_r({\bf k}_r)$, ${\bf k}_r\in[0,1]^{N-r}$. Then
for ${\bf k}_j\in[0,1]^{N-j}$ we define
\[\lb{014}
 \l_j({\bf k}_j)\ev\{\l:\ F_j(\l,{\bf k}_j)=0\ for\ \l\in\R\sm\bigcup_{r=0}^{j-1}\l_r([0,1]^{j-r},{\bf
 k}_j)\}.
\]
The condition after {\rm for} in \er{014} is natural because
otherwise some $F_r$ is zero and $F_{r+1}$ is not well defined (see
Defnition \ref{D1}). {\it We will call $\l_j$ Floquet-Bloch
dispersion branches for the operator $\cA$.} The following result
fully describes the spectrum of the operator $\cA$ and it also
allows us to recover the operator $\cA$ from its spectrum.

\begin{theorem}\lb{T3} i) Let $\cA$ be of the form \er{004}-\er{005}
with real scalar components $A_j$. Then any set $\l_j({\bf k}_j)$
\er{014} consists of no more than one element. The spectrum of $\cA$
is a union of intervals $\s_j$, $j<N$ and one eigenvalue $\s_N$
(some of $\s_j$ or $\s_N$ can be empty)
\[\lb{015}
 \s(\cA)=\bigcup_{r=0}^N\s_r,\ \ where\ \ \s_j=\l_j([0,1]^{N-j}).
\]
ii) Suppose that we have $N+1$ real continuous functions
$\l_r\ev\l_r({\bf k}_r)$, ${\bf k}_r\in[0,1]^{N-r}$, $r=0,...,N$
(this means that any set $\l_r({\bf k}_r)$ consists of one element)
satisfying
\[\lb{016}
 \l_j({\bf k}_j)\not\in\bigcup_{r=0}^{j-1}\l_r([0,1]^{j-r},{\bf
 k}_j),\ \ \forall j\ge1,\ \forall{\bf k}_j\in[0,1]^{N-j}.
\]
Then there exists the unique operator $\cA$ \er{004}-\er{005} with
Floquet-Bloch branches $\l_j$. Its components can be calculated by
induction $A_0\ev\l_0$ and $A_j({\bf k}_j)=-\langle
F_{j-1}^{-1}(\l_j({\bf k}_j),k_{j},{\bf k}_j)\rangle^{-1}_j$,
$j\ge1$. The spectrum of $\cA$ satisfies \er{015}.
\end{theorem}

{\bf Remark.} By analogy with \er{014}, the Floquet-Bloch dispersion
branches $\l_j({\bf k}_j)$ can also be defined for $M>1$ as zeroes
of $\det{\bf G}_j$, see \er{006}. The sets $\l_j({\bf k}_j)$ and
their projections $\l_r([0,1]^{j-r},{\bf
 k}_j)$ play important role in physics of waves, see,
e.g., \cite{K1}. The well-known "band structures" or band-gap
diagrams consist of the Floquet-Bloch branches. For different $j$
the components $\s_j=\l_j([0,1]^{N-j})$ have different physical (and
mathematical) nature. Roughly speaking, $\s_j$ characterizes the
waves that propagate along the defect of the dimension $N-j$ and
exponentially decrease in the perpendicular directions. The
corresponding waves are called guided and localized waves, see, e.g.
\cite{KKSP1}.

{\bf Example.} Let $M=1$, $N=2$. Suppose that we want to construct
the operator $\cA$ \er{004}-\er{005} with the spectrum
$\s=\cup\s_j$, where $\s_0=[0,1]$, $\s_1=[0.5,1.5]$, and $\s_2=2$.
In order to do this we can choose any real continuous functions
$\l_j$ satisfying \er{016} with $\l_j([0,1]^{N-j})=\s_j$, and use
after the procedure from Theorem \ref{T3}.ii) along with Definition
\ref{D1}. Ok, suppose that we want $\l_0=k_1k_2$, $\l_1=0.5+k_2$,
and $\l_2=2$. Then
$$
 A_0=k_1k_2,\ \
 A_1=-\langle
 F_0^{-1}(\l_1,k_1,k_2)\rangle^{-1}_1=-\lt(\int_0^1\frac{dk_1}{k_1k_2-0.5-k_2}\rt)^{-1}=\frac{k_2}{\ln(1+2k_2)},
$$
$$
 F_1(\l_2,k_2)=\frac{k_2}{\ln(1+2k_2)}+\lt(\int_0^1\frac{dk_1}{k_1k_2-2}\rt)^{-1}=\frac{k_2}{\ln(1+2k_2)}+\frac{k_2}{\ln(1-0.5k_2)},
$$
$$
 A_3=-\langle
 F_1^{-1}(\l_2,k_2)\rangle_2^{-1}=-\lt(\int_0^1\frac{\ln(1+2k_2)\ln(1-0.5k_2)dk_2}{k_2\ln(1+1.5k_2-k_2^2)}\rt)^{-1}=0.935...,
$$
$$
and\ then\
\cA=k_1k_2\cdot+\frac{k_2}{\ln(1+2k_2)}\langle\cdot\rangle_1+0.935...\langle\cdot\rangle_{1,2}.
$$
As a conclusion note that there are only few papers devoted to
applications of (finite or infinite) standard matrix-valued
continued fractions (MCF)  to spectral problems: in \cite{AP} some
general relations between Hamiltonians and MCF are presented; in
\cite{M} MCF are applied for calculating Green functions related to
some Hamiltonians; in \cite{SWP} the authors use MCF in analysis of
non-linear spectral problems; in \cite{O} some methods of obtaining
Floquet eigenvalues and eigensolutions based on MCF are discussed;
in \cite{L}, \cite{RK} the stability of the methods is analyzed.
Some applications of MCF to explicit representations of resolvent
operators can be found in \cite{I}, \cite{SI}. Scalar continued
fractions (CF) are used in Krein's inverse spectral problems, see,
e.g., \cite{T}. Also note an interesting connections between
orthogonal polynomials (including matrix-valued polynomials), CF,
and inverse spectral theory, see, e.g., \cite{GT}, \cite{DPS}.
Classical integral CF are introduced in \cite{MS}, \cite{BS} as
solutions of differential equations. They are also discussed in the
context of interpolation theory, see, e.g., \cite{M1}, \cite{MKM}.
At the same time, probably there are no papers devoted to integral
CF or MCF of the form \er{e010}.

The work is organized as follows. Section \ref{S1} contains the
proofs of Theorems \ref{T1}, \ref{T2}, and explicit resolvent
formulas. Section \ref{S2} contains the proof of Theorem \ref{T3}
Section \ref{S3} provides an application of our results to the
spectral problem of discrete Schr\"odinger operator acting on the
graphene with line and local inclusions. The conclusion is given in
Section \ref{S4}.

\section{Proofs of Theorems \ref{T1}, \ref{T2}, and explicit resolvent formulas}
\lb{S1}

The proof is based on the following result from \cite{Kjmaa}:

{\it Let $\cA$ be an operator of the form \er{001}-\er{003}. Then
the spectrum of $\cA$ is
\[\lb{100}
 \s(\cA)=\bigcup_{j=0}^N\s_j\ \ with\ \ \s_j=\{\l:\ \det{\bf E}_j=0\ for\ some\ {\bf
 k}_j\in[0,1]^{N-j}\},
\]
where matrix-valued functions ${\bf E}_j$ are defined as:
\[\lb{101}
 {\bf C}_0=({\bf A}_0-\l{\bf I})^{-1},\ \ {\bf E}_0={\bf C}_0^{-1},\ \ {\bf C}_1=-{\bf C}_0{\bf A}_1,\ \
 {\bf E}_1={\bf I}-\langle{\bf B}_1{\bf
 C}_1\rangle_{1},\ \ {\bf D}_1={\bf E}_1^{-1}{\bf B}_1{\bf C}_0
\]
and for $1<j\le N$ by induction
\[\lb{102}
 {\bf C}_j=-{\bf C}_0{\bf A}_j-\sum_{r=1}^{j-1}{\bf C}_r\langle{\bf D}_r{\bf
 A}_j\rangle_{1,r},\ \ \ {\bf E}_j={\bf I}-\langle{\bf B}_j{\bf
 C}_j\rangle_{1,j},
\]
\[\lb{103}
 {\bf D}_j={\bf E}_j^{-1}\lt({\bf B}_j{\bf C}_0+\sum_{r=1}^{j-1}\langle{\bf B}_j{\bf C}_{r}\rangle_{1,r}{\bf
 D}_{r}\rt).
\]}
Adapting \er{101}-\er{103} to the special case of $\cA$
\er{004}-\er{005} (where all ${\bf B}_r={\bf I}$) we obtain that:
\[\lb{104}
 {\bf E}_j={\bf I}-\langle{\bf
 C}_j\rangle_{1,j},\ \ {\bf D}_j={\bf E}_j^{-1}\lt({\bf C}_0+\sum_{r=1}^{j-1}\langle{\bf C}_{r}\rangle_{1,r}{\bf
 D}_{r}\rt),\ \ j\ge1,
\]
where $\sum_{r=1}^{0}$ is assumed to be zero (for $j=1$). Due to the
fact that ${\bf A}_r$ do not depend on $k_1,...,k_r$ we also deduce
\[\lb{105}
 {\bf C}_j=-{\bf H}_j{\bf
 A}_j,\ \ where\ \ {\bf H}_j={\bf C}_0+\sum_{r=1}^{j-1}{\bf C}_r\langle{\bf D}_r\rangle_{1,r},\ \ j\ge1.
\]
Using direct calculations we obtain explicit formulas for the first
few matrices ${\bf C}_0={\bf F}_0^{-1}$,  ${\bf E}_0={\bf F}_0$ and
\[\lb{106}
 {\bf H}_1={\bf
 F}_0^{-1},\ \ {\bf E}_1={\bf I}+\langle{\bf F}_0^{-1}\rangle_1{\bf
 A}_1=\langle{\bf F}_0^{-1}\rangle_1{\bf F}_1,\ \ {\bf D}_1={\bf
 F}_1^{-1}\langle{\bf F}_0^{-1}\rangle_1^{-1}{\bf F}_0^{-1},
\]
\[\lb{106a}
 {\bf H}_2={\bf
 F}_0^{-1}\langle{\bf
 F}_0^{-1}\rangle_1^{-1}{\bf F}_1^{-1},\ \ {\bf E}_2=\langle{\bf F}_1^{-1}\rangle_2{\bf F}_2,\ \ {\bf D}_2={\bf
 F}_2^{-1}\langle{\bf F}_1^{-1}\rangle_2^{-1}{\bf F}_1^{-1}\langle{\bf F}_0^{-1}\rangle_1^{-1}{\bf F}_0^{-1},
\]
where the matrices ${\bf F}$ are defined in \er{e010}. Now, we prove
by induction the following identities
\[\lb{107}
 {\bf E}_j=\langle{\bf
 F}_{j-1}^{-1}\rangle_j{\bf F}_j,\ \ {\bf D}_j={\bf
 F}_j^{-1}\lt(\langle{\bf F}_{j-1}^{-1}\rangle_j^{-1}{\bf
 F}_{j-1}^{-1}\rt)...\lt(\langle{\bf F}_{0}^{-1}\rangle_1^{-1}{\bf
 F}_{0}^{-1}\rt)
\]
and
\[\lb{108}
 {\bf H}_j=\lt({\bf
 F}_{0}^{-1}\langle{\bf F}_{0}^{-1}\rangle_1^{-1}\rt)...\lt({\bf
 F}_{j-2}^{-1}\langle{\bf F}_{j-2}^{-1}\rangle_{j-1}^{-1}\rt){\bf
 F}_{j-1}^{-1}.
\]
Really, \er{105} gives us
\[\lb{109}
 {\bf H}_{j+1}={\bf H}_j+{\bf C}_j\langle{\bf D}_j\rangle_{1,j}={\bf
 H}_j({\bf I}-{\bf A}_j\langle{\bf D}_j\rangle_{1,j}).
\]
Using \er{107}, Fubini's theorem, and the fact that ${\bf F}_r$ do
not depend on $k_1,...,k_r$ (see \er{e010}) we obtain
\[\lb{110}
 \langle{\bf D}_j\rangle_{1,j}={\bf
 F}_j^{-1}\lt\langle\langle{\bf F}_{j-1}^{-1}\rangle_j^{-1}{\bf
 F}_{j-1}^{-1}...\lt\langle\langle{\bf F}_{0}^{-1}\rangle_1^{-1}{\bf
 F}_{0}^{-1}\rt\rangle_1\rt\rangle_{2,j}
\]
\[\lb{111}
 ={\bf
 F}_j^{-1}\lt\langle\langle{\bf F}_{j-1}^{-1}\rangle_j^{-1}{\bf
 F}_{j-1}^{-1}...\langle{\bf F}_{0}^{-1}\rangle_1^{-1}\langle{\bf
 F}_{0}^{-1}\rangle_1\rt\rangle_{2,j}={\bf
 F}_j^{-1}\lt\langle\langle{\bf F}_{j-1}^{-1}\rangle_j^{-1}{\bf
 F}_{j-1}^{-1}...\rt\rangle_{2,j}=...={\bf
 F}_j^{-1}.
\]
Substituting \er{110}-\er{111} into \er{109} and using Definition
\ref{D1} we get
\[\lb{112}
 {\bf H}_{j+1}={\bf
 H}_j({\bf I}-{\bf A}_j{\bf
 F}_j^{-1})={\bf
 H}_j\langle{\bf F}_{j-1}^{-1}\rangle_j^{-1}{\bf
 F}_j^{-1}
\]
which confirms \er{108} for $j+1$. The similar arguments applied to
${\bf E}_{j+1}$ (see \er{104}) give us
\[\lb{113}
 {\bf E}_{j+1}={\bf I}-\langle{\bf
 C}_{j+1}\rangle_{1,j+1}={\bf I}+\langle{\bf
 H}_{j+1}\rangle_{1,j+1}{\bf A}_{j+1}={\bf I}+\langle{\bf
 F}_j^{-1}\rangle_{j+1}{\bf A}_{j+1}=\langle{\bf
 F}_j^{-1}\rangle_{j+1}{\bf F}_{j+1}
\]
which confirms the first identity in \er{107} for $j+1$. Identities
\er{104} lead to
\[\lb{114}
 {\bf E}_{j+1}{\bf D}_{j+1}={\bf E}_j{\bf D}_j+\langle{\bf
 C}_{j}\rangle_{1,j}{\bf D}_{j}={\bf D}_j.
\]
Thus, by \er{107} (it is already proved for ${\bf E}_{j+1}$, see
\er{113}) we have
\[\lb{115}
 {\bf D}_{j+1}={\bf E}_{j+1}^{-1}{\bf D}_j={\bf F}_{j+1}^{-1}\langle{\bf
 F}_j^{-1}\rangle_{j+1}^{-1}{\bf D}_j
\]
which confirms the second identity in \er{107} for $j+1$. Now,
identities \er{107}-\er{108} are completely proved for all
$j=2,...,N$. We obtain also that ${\bf G}_j={\bf E}_j$ (see
\er{107}, \er{007}). Due to \er{100} this means that Theorem
\ref{T1} is true. Identities \er{009} are obvious by the definition
of ${\bf F}$ \er{e010}. Let us compute the second identity in
\er{008}
\[\lb{116}
 \left\langle({\bf G}_0...{\bf
 G}_{j-1})^{-1}\right\rangle_{1,j}^{-1}({\bf G}_j-{\bf I})=\left\langle({\bf E}_0...{\bf
 E}_{j-1})^{-1}\right\rangle_{1,j}^{-1}({\bf E}_j-{\bf I})=
\]
\[\lb{117}
 \langle{\bf D}_{j-1}\rangle_{1,j}^{-1}(\langle{\bf
 F}_{j-1}^{-1}\rangle_j{\bf F}_j-{\bf I})=\langle{\bf
 F}_{j-1}^{-1}\rangle_j^{-1}(\langle{\bf
 F}_{j-1}^{-1}\rangle_j{\bf F}_j-{\bf I})={\bf A}_j,
\]
where we have used \er{107}, \er{110}-\er{111}, and Definition
\ref{D1}. Thus, Theorem \ref{T2} is proved.

Identities \er{011},\er{012} follow from \er{105},\er{112},\er{115},
and from the general result about the resolvent (see \cite{Kjmaa})
\[\lb{118}
 (\cA-\l\cI)^{-1}={\bf C}_0\cdot+{\bf C}_1\langle{\bf
 D}_1\cdot\rangle_1+...+{\bf C}_N\langle{\bf D}_N\cdot\rangle_{1,N}.
\]
If all ${\bf A}_j$ are self-adjoint then all ${\bf F}_j$ are
self-adjoint and then ${\bf D}_j^*={\bf H}_{j+1}$, see
\er{112},\er{115}. These arguments lead to \er{013}.

\section{Proof of Theorem \ref{T3}}

\lb{S2}

We can consider $\l\in\R$ only because the operator $\cA$ is
self-adjoint. For any fixed ${\bf k}^0\in[0,1]^N$ the function
$F_0(\l,{\bf k}_0)=A_0({\bf k}_0)-\l$ is an injection on $\l\in\R$.
For any fixed ${\bf k}_1\in[0,1]^{N-1}$ the function $\langle
F_0^{-1}(\l,k_1,{\bf k}_1)\rangle_1^{-1}$ does not have zeroes in
the set $\l\in\R\sm\l_0([0,1],{\bf k}_1)=\R\sm A_0([0,1],{\bf k}_1)$
because $F^{-1}_0(\l,k_1,{\bf k}_1)$ does not change sign inside
this set. Thus, by Theorem \ref{T1} the spectral component
$$
 \s_1=\{\l:\ F_1(\l,{\bf k}_1)=0\ for\ some\ {\bf
 k}_1\in[0,1]^{N-1}\}\ev\l_1([0,1]^{N-1}).
$$
Consider some fixed ${\bf k}_1\in[0,1]^{N-1}$ and two different
points $x_1\ne x_2$ belonging to the set $\R\sm\l_0([0,1],{\bf
k}_1)$. Then $F_0(x_1,k_1,{\bf k}_1)\ne F_0(x_2,k_1,{\bf k}_1)$ and
$$
 F_0(x_1,k_1,{\bf k}_1)<(or\ >) F_0(x_2,k_1,{\bf k}_1)\ \Rightarrow\
 \langle F_0^{-1}(x_1,k_1,{\bf k}_1)\rangle_1^{-1}<(or\ >)\langle F_0^{-1}(x_2,k_1,{\bf k}_1)\rangle_1^{-1}\ \Rightarrow\
$$
$$
 F_1(x_1,{\bf k}_1)<(or\ >) F_1(x_2,{\bf k}_1),
$$
where we use the fact that  for any fixed $k_1\in[0,1]$, ${\bf
k}_1\in[0,1]^{N-1}$ the continuous function $F_0(\l,k_1,{\bf k}_1)$
is an injection on the set $\l\in\R\sm\l_0([0,1],{\bf k}_1)$ and
hence $F_0(x_1,k_1,{\bf k}_1)-F_0(x_2,k_1,{\bf k}_1)$ is not zero
and has the same sign for all $k_1\in[0,1]$. Thus, for any fixed
${\bf k}_1\in[0,1]^{N-1}$ the function $F_1(\l,{\bf k}_1)$ is an
injection on the set $\l\in\R\sm\l_0([0,1],{\bf k}_1)$. Then the set
$$
 \l_1({\bf k}_1)=\{\l:\ F_1(\l,{\bf k}_1)=0\}
$$
consists of no more than one point. The same arguments allow us to
finish the proof of i) by induction.

The results of ii) follow also from the facts: for any fixed ${\bf
k}_j$ the function $F_j(\l,{\bf k}_j)$ is an injection (it can be
proved by induction, see above) on the set
$\l\in\R\sm\bigcup_{r=0}^{j-1}\l_r([0,1]^{j-r},{\bf
 k}_j)\}$, and
$$
 F_j(\l,{\bf k}_j)=A_j({\bf k}_j)+\langle F^{-1}_{j-1}(\l,k_{j},{\bf
 k}_j)\rangle_j^{-1}=0
$$
for $\l=\l_j({\bf k}_j)$ if and only if $A_j({\bf k}_j)=-\langle
F^{-1}_{j-1}(\l,k_{j},{\bf
 k}_j)\rangle_j^{-1}$.

\section{Example} \lb{S3}

\begin{figure}[h]
\center{\includegraphics[width=0.65\linewidth]{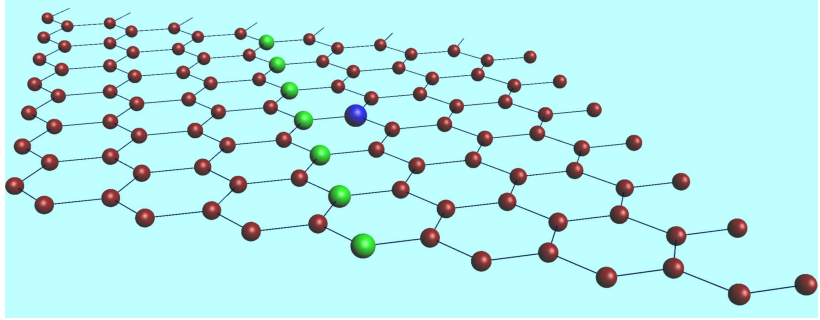}}
\caption{Graphene with 1D line defect and one local defect.}
\label{fig2}
\end{figure}

\begin{figure}[h]
\begin{minipage}[h]{0.49\linewidth}
\center{\includegraphics[width=0.9\linewidth]{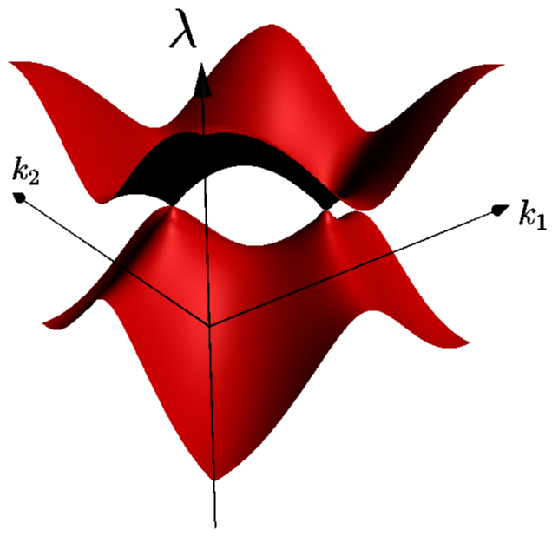}} \\ (a)
\end{minipage}
\hfill
\begin{minipage}[h]{0.49\linewidth}
\center{\includegraphics[width=0.9\linewidth]{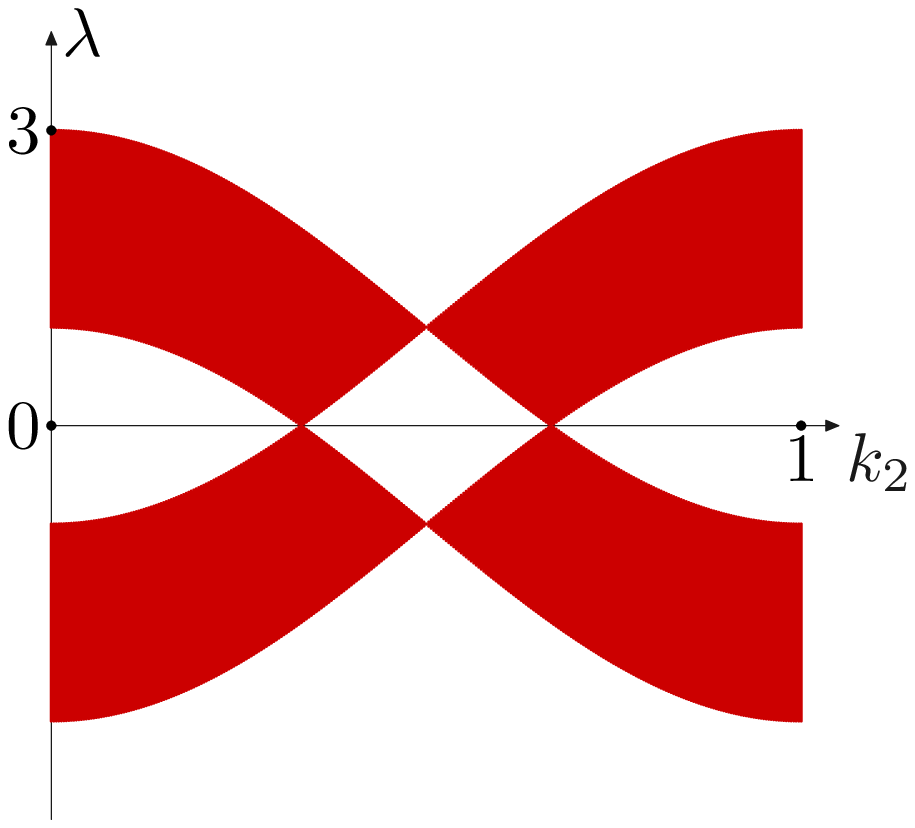}} \\ (b)
\end{minipage}
\caption{(a) Propagating dispersion surfaces $\l_{\pm}$ \er{307} and
(b) their projections (defined by the curves $\l_{\pm\min,\pm\max}$
\er{309}).} \label{fig3}
\end{figure}

\begin{figure}[h]
\begin{minipage}[h]{0.49\linewidth}
\center{\includegraphics[width=0.9\linewidth]{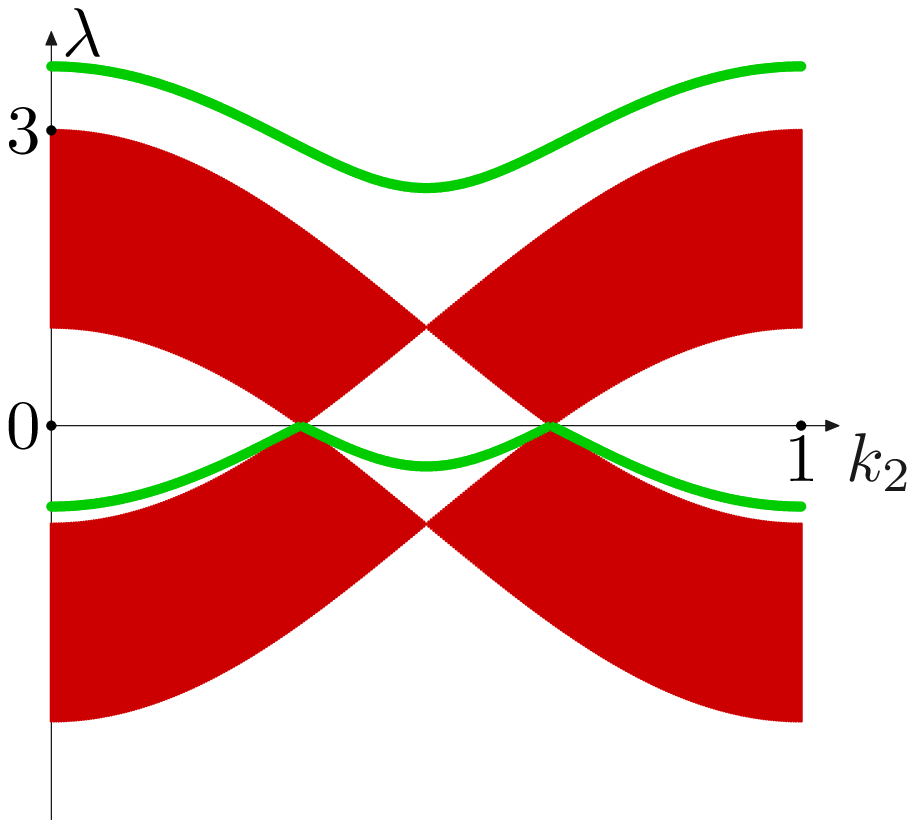}} \\
(a)
\end{minipage}
\hfill
\begin{minipage}[h]{0.49\linewidth}
\center{\includegraphics[width=0.9\linewidth]{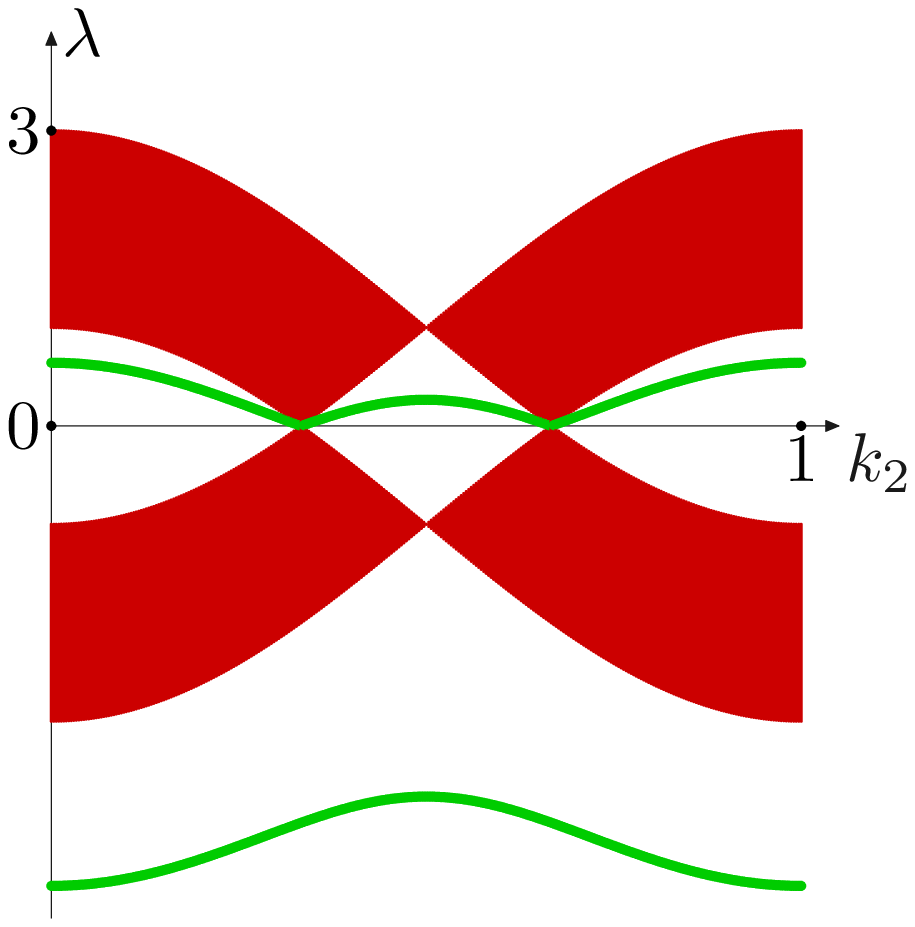}} \\
(b)
\end{minipage}
\caption{Projections of propagating dispersion surfaces $\l_{\pm}$
\er{307} (the same as in Fig. \ref{fig3}.(b)) and guided dispersion
curves $\l(k_2)$ \er{310} for the potentials (a) $V_1=2$, (b)
$V_1=-3.5$.} \label{fig4}
\end{figure}

\begin{figure}[h]
\begin{minipage}[h]{0.49\linewidth}
\center{\includegraphics[width=0.9\linewidth]{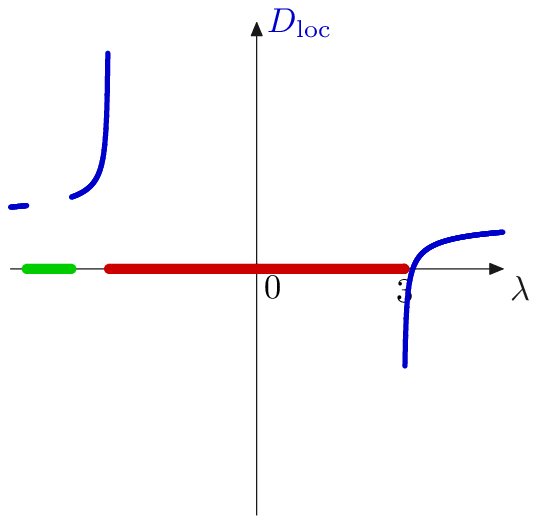}} \\
(a)
\end{minipage}
\hfill
\begin{minipage}[h]{0.49\linewidth}
\center{\includegraphics[width=0.9\linewidth]{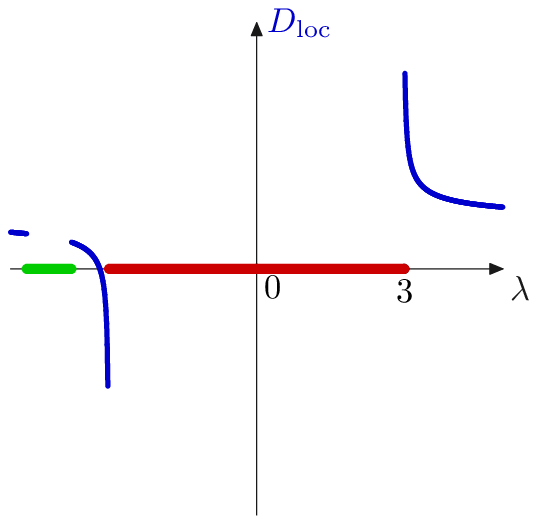}} \\
(b)
\end{minipage}
\caption{Propagating $\s_0$ (red line), guided $\s_1$ (green line)
spectral components, and the function $D_{\rm loc}(\l)$ \er{311}
computed for the potential $V_1=-3.5$ (see Fig. \ref{fig4}.(b)) and
the potentials (a) $V_2=200$, (b) $V_2=-200$. Zeroes of $D_{\rm
loc}(\l)$ are eigenvalues located in the gaps of the continuous
spectrum, i.e. in $\R\sm(\s_0\cup\s_1)$.} \label{fig5}
\end{figure}

In this example we study the spectral problem for the discrete
Schr\"odinger operator acting on the graphene with a line defect and
a local defect, see Fig. \ref{fig2}. The graphene $\G$ is
$2$-periodic lattice with $2$-point unit cell, the line defect
$\G_1$ is $1$-periodic sublattice, and the local defect is one
point, i.e.
\[\lb{300}
 \G=\{1,2\}\ts\Z^2,\ \ \G_1=\{1\}\ts0\ts\Z,\ \ \G_2=\{2\}\ts0\ts0.
\]
Consider the discrete Schr\"odinger operator $\cA$ acting on the
graphene
\[\lb{301}
 \cA:\ell^2(\G)\to\ell^2(\G),\ \ \cA U_{\bf n}=\wt\D U_{\bf n}+V_{\bf n}U_{\bf n},\ \ {\bf
 n}\in\G,
\]
where the modified Laplace operator $\wt\D U_{\bf n}=\sum_{{\bf
n}'\sim{\bf n}}U_{{\bf n}'}$ is the sum of the values of the wave
function $U\in\ell^2(\G)$ counting at neighbor points (for each
point ${\bf n}$ there are $3$ neighbor points ${\bf n}'$). The
potential is
\[\lb{302}
 V_{\bf n}=\ca 0, & {\bf n}\in\G\sm(\G_1\cup\G_2),\\ V_1, & {\bf n}\in\G_1,\\ V_2, & {\bf n}\in\G_2.\ac
\]
After applying FFB transformation
\[\lb{303}
 \cF:\ell^2(\G)\to L^2_{2,2},\ \
 \cF(U_{\bf n})=\sum_{{\bf m}\in\Z^2} e^{2\pi i{\bf m}^*{\bf k}}\ma U_{1\ts{\bf m}} \\ U_{2\ts{\bf m}}
 \am,\ \ {\bf k}\in[0,1]^2
\]
the Schr\"odinger operator \er{301} takes the form
\[\lb{304}
 \hat\cA=\cF\cA\cF^{-1}:L^2_{2,2}\to L^2_{2,2},\ \ \hat\cA={\bf
 A}_0\cdot+{\bf A}_1\langle\cdot\rangle_1+{\bf
 A}_2\langle\cdot\rangle_{1,2},
\]
where
\[\lb{305}
 {\bf A}_0=\ma 0 & e^{-2\pi i k_1}(1+e^{2\pi i k_2})+1 \\
   e^{2\pi i k_1}(1+e^{-2\pi i k_2})+1 & 0\am,
\]
\[\lb{305a}
   {\bf A}_1=\ma V_1 & 0 \\ 0 & 0 \am,\ \ {\bf A}_2=\ma 0 & 0 \\ 0 & V_2 \am,
\]
and ${\bf I}$ is $2\ts2$ identity matrix. The spectral component
$\s_0$ (see \er{006}-\er{007}) is determined by
\[\lb{306}
 \det{\bf G}_0=\det{\bf F}_0=\det({\bf A}_0-\l{\bf I})=0
\]
which leads to two Floquet-Bloch dispersion branches (surfaces)
\[\lb{307}
 \l=\l_{\pm}=\pm\sqrt{3+2\cos2\pi k_1+2\cos2\pi k_2+2\cos2\pi
 (k_1-k_2)}.
\]
Then the spectral component $\s_0=\l_-([0,1]^2)\cup\l_+([0,1]^2)$.
The surfaces $\l_{\pm}({\bf k})$, ${\bf k}=(k_1,k_2)\in[0,1]^2$ are
plotted in Fig. \ref{fig3}. Up to elementary coordinate ($k_1\to
k_1,k_2\to k_1-k_2$ which does not have affect on the spectrum)
transformations, they are the same as in \cite{KP}. Each point of
the surfaces corresponds to the bounded non-attenuated wave-function
$U_{\bf n}$ (so-called propagating eigenfunction). That is why these
dispersion surfaces are called propagating dispersion surfaces.
While the propagating spectrum is usually easy to compute other
components are more complex. The spectral component $\s_1$ (see
\er{006}-\er{007}) is determined by the determinant of the matrix
\[\lb{308}
 {\bf G}_1={\bf I}+\langle({\bf A}_0-\l{\bf I})^{-1}\rangle_1{\bf
 A}_1=\ma 1\pm\frac{V_1\l}{\sqrt{(\l^2-1-4\cos^2\pi k_2)^2-16\cos^2\pi k_2}} & ... \\ 0 &
 1\am,
\]
where the sign $\pm$ depends on whether $\l$ is above or below the
projection of the propagating branches $\l_{\pm}({\bf k})$ on the
plane $(\l,k_2)$. This projection consists of areas bounded by the
four curves
\[\lb{309}
 \l_{\pm\min,\pm\max}(k_2)=\pm\sqrt{3+2\cos2\pi k_2\pm4\cos\pi k_2}.
\]
Returning to the spectral component $\s_1$, we have that it is
determined by the condition $\det{\bf G}_1=0$ (see \er{006}-\er{007}
and \er{308}) which leads to
\[\lb{310}
 \l^2=\l^2(k_2)=\frac{2(1+4\cos^2\pi k_2)+V_1^2\pm\sqrt{64\cos^2\pi k_2+4(1+4\cos^2\pi
 k_2)V_1^2+V_1^4}}2,
\]
where the sign $\pm$ and the sign of $\l(k_2)$ depend on the
location of $\l(k_2)$ (below or above the curves
$\l_{\pm\min,\pm\max}(k_2)$) and on the sign of the potential $V_1$.
The spectral component $\s_1$ is now $\s_1=\l([0,1])$. Each point
$\l(k_2)$ (where $k_2\in[0,1]$) corresponds to the bounded
wave-function $U_{\bf n}$ which is non-attenuated along the line
defect and exponentially decreasing in the perpendicular directions
(so-called guided eigenfunction). That is why these dispersion
curves are called guided dispersion curves. They are plotted in Fig.
\ref{fig4}.
For the discrete spectrum $\s_2$ we need to calculate zeroes of the
following function (see \er{006}-\er{007})
\[\lb{311}
 D_{\rm loc}(\l)=\det{\bf G}_2=\det\lt({\bf I}+\lt\langle\lt({\bf
 A}_1+\left\langle\frac{\bf I}{{\bf A}_0-\l{\bf
 I}}\right\rangle_1^{-1}\rt)^{-1}\rt\rangle_2{\bf A}_2\rt).
\]
This function can not be expressed in elementary functions but it
can easily be calculated numerically. On the Fig. \ref{fig5}, there
are some examples of presence of eigenvalues.

\section{Conclusion} \lb{S4}

We have shown that the set of matrix-valued integral continued
fractions defined by the components of the periodic operator with
defects determine the spectrum, the resolvent, and, of course, the
components of the operator explicitly. Roughly speaking, this set is
enough to solve the various direct and inverse spectral problems.

\section*{Acknowledgements}
This work was supported by the RSF project
N\textsuperscript{\underline{o}}15-11-30007. I would also like to
thank Dr. Royer for useful discussions.

\bibliography{contfrac}

\end{document}